\documentclass[12pt]{amsart}
\usepackage{amssymb,amsmath,amsfonts,latexsym,setspace}
\usepackage{bm}
\newcommand{\newsection}[1]{\setcounter{equation}{0} \section{#1}}
\newcommand{\bea}{\begin{eqnarray}}
\newcommand{\eea}{\end{eqnarray}}

\newcommand{\clb}{\mathcal{B}}

\newcommand{\cld}{\mathcal{D}}
\newcommand{\cle}{\mathcal{E}}
\newcommand{\clf}{\mathcal{F}}

\newcommand{\clh}{\mathcal{H}}
\newcommand{\clk}{\mathcal{K}}

\newcommand{\clq}{\mathcal{Q}}

\newcommand{\cls}{\mathcal{S}}

\newcommand{\raro}{\rightarrow}

\def \qed {\hfill \vrule height6pt width 6pt depth 0pt}
\def\textmatrix#1&#2\\#3&#4\\{\bigl({#1 \atop #3}\ {#2 \atop #4}\bigr)}
\def\dispmatrix#1&#2\\#3&#4\\{\left({#1 \atop #3}\ {#2 \atop #4}\right)}
\newcommand{\be}{\begin{equation}}
\newcommand{\ee}{\end{equation}}
\newcommand{\ben}{\begin{eqnarray*}}
\newcommand{\een}{\end{eqnarray*}}

\newcommand{\NI}{\noindent}

\newcommand{\bi}{\begin{itemize}}
\newcommand{\ei}{\end{itemize}}

\newtheorem{Theorem}{\sc Theorem}[section]
\newtheorem{Lemma}[Theorem]{\sc Lemma}
\newtheorem{Proposition}[Theorem]{\sc Proposition}
\newtheorem{Corollary}[Theorem]{\sc Corollary}
\newtheorem{Definition}[Theorem]{\sc Definition}
\newtheorem{Example}[Theorem]{\sc Example}
\newtheorem{Remark}[Theorem]{\sc Remark}
\newtheorem{Note}[Theorem]{\sc Note}
\newtheorem{Question}{\sc Question}
\newtheorem{ass}[Theorem]{\sc Assumption}
\newcommand{\bt}{\begin{Theorem}}
\def\beginlem{\begin{Lemma}}
\def\beginprop{\begin{Proposition}}
\def\begincor{\begin{Corollary}}
\def\begindef{\begin{Definition}}
\def\beginexamp{\begin{Example}}
\def\beginrem{\begin{Remark}}
\def\beginq{\begin{Question}}
\def\beginass{\begin{ass}}
\def\beginnote{\begin{Note}}
\newcommand{\et}{\end{Theorem}}
\def\endlem{\end{Lemma}}
\def\endprop{\end{Proposition}}
\def\endcor{\end{Corollary}}
\def\enddef{\end{Definition}}
\def\endexamp{\end{Example}}
\def\endrem{\end{Remark}}
\def\endq{\end{Question}}
\def\endass{\end{ass}}
\def\endnote{\end{Note}}

\begin{document}

\title[Operator Theory on Symmetrized Bidisc]{Operator Theory on
Symmetrized Bidisc}

\author[Jaydeb Sarkar]{Jaydeb Sarkar}
\address{Indian Statistical Institute, Statistics and Mathematics Unit, 8th Mile, Mysore Road, Bangalore, 560059, India}
\email{jay@isibang.ac.in, jaydeb@gmail.com}




\begin{abstract}
A commuting pair of operators $(S, P)$ on a Hilbert space $\clh$ is
said to be a $\Gamma$-contraction if the symmetrized bidisc
\[\Gamma = \{(z_1+z_2, z_1 z_2) :|z_1|, |z_2| \leq 1\}\]is a
spectral set of the tuple $(S, P)$. In this paper we develop some
operator theory inspired by Agler and Young's results on a model
theory for $\Gamma$-contractions.

\NI We prove a Beurling-Lax-Halmos type theorem for
$\Gamma$-isometries. Along the way we solve a problem in the
classical one-variable operator theory, namely, a non-zero
$M_z$-invariant subspace $\cls$ of $H^2_{\cle_*}(\mathbb{D})$ is
invariant under the analytic Toeplitz operator with the
operator-valued polynomial symbol $p(z) = A + A^* z$ if and only if
the Beurling-Lax-Halmos inner multiplier $\Theta$ of $\cls$
satisfies
\[(A + A^*z) \Theta = \Theta (B + B^*z),\]for some unique operator
$B$.

\NI We use a "pull back" technique to prove that a completely
non-unitary $\Gamma$-contraction $(S, P)$ can be dilated to a pair
\[(((A + A^*M_z) \oplus U), (M_z \oplus M_{e^{it}})),\]which is the
direct sum of a $\Gamma$-isometry and a $\Gamma$-unitary on the
Sz.-Nagy and Foias functional model of $P$, and that $(S, P)$ can be
realized as a compression of the above pair in the functional model
$\clq_P$ of $P$ as
\[(\bm{P}_{\clq_P} ((A+A^*M_z) \oplus U)|_{\clq_P}, \bm{P}_{\clq_P} (M_z \oplus
M_{e^{it}})|_{\clq_P}).\] Moreover, we show that this representation
is unique. We prove that a commuting tuple $(S, P)$ with $\|S\| \leq
2$ and $\|P\| \leq 1$ is a $\Gamma$-contraction if and only if there
exists a compressed scalar operator $X$ with the decompressed
numerical radius not greater than one such that
\[S = X + P X^*.\] In the commutant lifting set up, we obtain a unique and explicit
solution to the lifting of $S$ where $(S, P)$ is a completely
non-unitary $\Gamma$-contraction. Our results concerning the
Beurling-Lax-Halmos theorem of $\Gamma$-isometries and the
functional model of $\Gamma$-contractions answers a pair of
questions of J. Agler and N. J. Young.
\end{abstract}

\maketitle

\newsection{Introduction}

The notion of spectral set was introduced by J. von Neumann in
\cite{vN} where he proved that the closed unit disk
$\overline{\mathbb{D}}$, where $\mathbb{D} = \{z \in \mathbb{C}: |z|
< 1\}$, is a spectral set of a bounded linear operator on a Hilbert
space if and only if the operator is a contraction. Later in
\cite{Nagy}, Sz.-Nagy proved that a bounded linear operator is a
contraction if and only if the operator has a unitary dilation.
Therefore von Neumann's result can be derived from Sz.-Nagy's
unitary dilation. Since then, one of the most celebrated problems in
operator theory is to determine the class of commuting $n$-tuple of
operators for which a normal $\partial K$-dilation exists, where $K
\subseteq \mathbb{C}^n$ is compact and $n \geq 1$. We recall that a
commuting tuple $(T_1, \ldots, T_n)$ on $\clh$ has a normal
$\partial K$-dilation if there exists a tuple of commuting normal
operators $(N_1, \ldots, N_n)$ on $\clk \supseteq \clh$ such that
$\sigma_T(N_1, \ldots, N_n) \subseteq \partial K$ and
\[\bm{P}_{\clh} p(N_1, \ldots, N_n)|_{\clh} = p(T_1, \ldots, T_n),\] for
all $p \in \mathbb{C}[z_1, \ldots, z_n]$. Here $\bm{P}_{\clh}$ is
the orthogonal projection of $\clk$ onto $\clh$. Many studies in
this problem have been carried so far. In particular, it is known
that the normal $\partial K$-dilation holds if $K$ is the closure of
an annulus \cite{A} and fails when $K$ is a triply connected domain
in $\mathbb{C}$ \cite{DM}. The theory becomes more subtle when the
spectral set is assumed to be a subset of $\mathbb{C}^n$ ($n
> 1$).

On the other hand, it is well known that for $n \geq 2$, the von
Neumann's inequality fails in general with the exception that a pair
of commuting contractions can be dilated to a pair of commuting
unitary operators \cite{Ando}. One versions of von Neumann's
inequality for domains like ball and general symmetric domains
require to replace the sup norm of the polynomials by operator norm
of certain natural multiplier algebras. Now, we define the notion of
$\Gamma$-contractions.

A pair of commuting operators $(S, P)$ on a Hilbert space $\clh$ is
said to be a \textit{$\Gamma$-contraction} if the symmetrized bidisc
\[\Gamma = \{(z_1 + z_2, z+_1 z_2) : |z_1|, |z_2| \leq 1\}.\] is a
spectral set of $(S, P)$. That is, for all polynomial $p \in
\mathbb{C}[z_1, z_2]$, \[\|p(S, P)\| \leq \sup_{z \in \Gamma}
|p(z)|.\]

\NI In particular, if $(S, P)$ is a $\Gamma$-contraction then $\|S\|
\leq 2$ and $\|P\| \leq 1$. Note also that the symmetrized bidisc
$\Gamma$ is the range of $\pi$ restricted to the closed bidisc
$\bar{\mathbb{D}}^2$ where $\pi : \mathbb{C}^2 \raro \mathbb{C}^2$
is the proper holomorphic map defined by
\[\pi(z_1, z_2) = (z_1+z_2, z_1 z_2),\]for all $(z_1, z_2) \in
\mathbb{C}^2$.

There is a significant difference between $\Gamma$ and other bounded
symmetric domains considered earlier by many researches in the
development of analytic model theory (cf. \cite{AEM}). For instance,
$\Gamma$ is polynomially convex \cite{AY2} but non-convex and
inhomogeneous. This in turns makes the theory of
$\Gamma$-contraction more appealing and useful in the study of the
classical and several variables operator theory.

In \cite{AY5}, Agler and Young developed a $\Gamma$-isometric
dilation theory for $\Gamma$-contractions. In this paper, we develop
an explicit $\Gamma$-isometric dilation and functional model of
$\Gamma$-contractions. Furthermore, we provide a characterization of
$\Gamma$-contractions which is compatible with the geometry of the
domain $\Gamma$. Moreover, we obtain a characterization of invariant
subspaces of $\Gamma$-isometries.

Our results improve, generalize and unify some recent known results
(\cite{BPR}, \cite{BP}) and answer positively a pair of problems
mentioned by Agler and Young (page 58, \cite{AY5}).

Finally, it is worth mentioning that our classification results (cf.
Theorems \ref{XP} and \ref{geom-2}) imply that we must confine
ourselves to the class of contractions and compressed scalar
operators (see Section 6) in order to have a concrete realization of
operator tuples with $\Omega \subseteq \mathbb{C}^n$ as a spectral
set, where $\Omega$ is the symmetrized polydisc or a higher
dimensional domain related to the symmetrized bidisc (cf.
\cite{AWY}).

The organization of the paper is as follows.

In Section 2, we recall some basic definitions and results in the
theory of $\Gamma$-contractions.

In Section 3, we provide some basic classification results of pure
$\Gamma$-isometries. We obtain a Beurling-Lax-Halmos type theorem
characterizing joint invariant subspaces of a pure
$\Gamma$-isometry.

In Section 4, we prove a factorization result concerning isometric
dilation of a completely non-unitary contraction. Moreover, we use a
"pull back" argument to the factorization and obtain a functional
model for completely non-unitary $\Gamma$-contractions.

In Section 5, we show that the functional model of a completely non
unitary $\Gamma$-contraction is unique.

In Section 6, we proceed to a new characterization of
$\Gamma$-contractions.

In Section 7, we conclude with a number of results and remarks
concerning $\Gamma$-isometric Hardy modules, isometrically
isomorphic submodules of $\Gamma$-isometric Hardy modules and a
solution to the commutant lifting problem.

In this paper, all Hilbert spaces are assumed to be separable and
over the field of complex numbers.

\newsection{Preliminaries}
In this section, we will gather together some of the necessary
definitions and results on $\Gamma$-contractions which we will
employ later in the paper. For more details about
$\Gamma$-contractions, we refer readers to the seminal work of Agler
and Young \cite{AY1}-\cite{AY6} (also see \cite{BPR}, \cite{BP}).

In what follows, we shall denote a pair of commuting operators by
$(S, P)$, for ``sum'' and ``product''. However, it is far from true
that a $\Gamma$-contraction is necessarily a sum and product of a
pair of commuting contractions.

Let $(S, P)$  be a pair of commuting operators on a Hilbert space
$\clh$. Then $(S, P)$ is said to be

(i) \textit{$\Gamma$-unitary} if $S$ and $P$ are normal operators
and the joint spectrum $\sigma(S, P)$ is contained in the
distinguished boundary of $\Gamma$.

(ii) \textit{$\Gamma$-isometry} if $(S, P)$ has a $\Gamma$-unitary
extension.

(iii) \textit{$\Gamma$-co-isometry} if $(S^*, P^*)$ is a
$\Gamma$-isometry.

The following theorem is due to Agler and Young \cite{AY5}.

\begin{Theorem}\label{gamma-unitary}
\textsf{(Agler and Young)} Let $(S,P)$ be a pair of commuting
operators on a Hilbert space $\clh$. Then the following statements
are equivalent:

(i) $(S, P)$ is a $\Gamma$-unitary.

(ii) $P$ is unitary and $S = S^*P$ and $\|S\| \leq 2$.

(iii) There exists commuting unitary operators $U_1, U_2$ on $\clh$
such that \[S = U_1 + U_2, \quad P = U_1 U_2.\]
\end{Theorem}
Note that for a $\Gamma$-isometry $(S, P)$ on $\clh$ we have
\[S = \tilde{S}|_{\clh} \quad \mbox{and} \quad P =
\tilde{P}|_{\clh},\]where $(\tilde{S}, \tilde{P})$ is a
$\Gamma$-unitary on a Hilbert space $\clk \supseteq \clh$.
Consequently, a necessary condition for a pair of commuting
operators $(S, P)$ to be a $\Gamma$-isometry is that $P$ is an
isometry. A $\Gamma$-isometry $(S,P)$ is said to be \textit{pure
$\Gamma$-isometry} if $P$ is a pure isometry, that is, $P$ does not
have any unitary part. The following characterization result play an
important role in the sequel.

\begin{Theorem}\label{gamma-isometry}
\textsf{(Agler and Young)} Let $S, P$ be commuting operators on a
Hilbert space $\clh$. Then $(S,P)$ is a pure $\Gamma$-isometry if
and only if there exists a Hilbert space $\cle$, a unitary operator
$U : \clh \raro H^2_{\cle}(\mathbb{D})$ and $A \in \clb(\cle)$ such
that $w(A) \leq 1$ and \[S = U^* M_{\varphi} U, \quad P =
U^*M_zU,\]where
\[\varphi(z) = A + A^*z, \quad z \in \mathbb{D}.\]
\end{Theorem}

Here $w(A)$ is the numerical radius of the operator $A \in
\clb(\cle)$, that is,
\[w(A) = \sup \{|\langle Ah, h \rangle| : h \in \cle, \|h\| \leq 1\}.\]

A contraction $P$ on $\clh$ is said to be \textit{completely
non-unitary} (or \textit{c.n.u.}) if there is no non-zero
$P$-reducing subspace $\clh_u \subseteq \clh$ such that
$T|_{\clh_u}$ is unitary. It is known that a contraction $P$ on
$\clh$ can be uniquely decomposed as $P = P|_{\clh_n} \oplus
P|_{\clh_u}$ where $\clh_n$ and $\clh_u$ are $P$-reducing subspaces
of $\clh$ and $P|_{\clh_n}$ is a c.n.u. contraction and
$P|_{\clh_u}$ is a unitary contraction. Moreover, let $(S, P)$ be a
$\Gamma$-contraction for some operator $S$ on $\clh$. Then $\clh_n$
and $\clh_u$ are $S$-reducing too and $(S|_{\clh_u}, P|_{\clh_u})$
is a $\Gamma$-unitary and $(S|_{\clh_n}, P|_{\clh_n})$ is a
$\Gamma$-contraction (Theorem 2.8 in \cite{AY5}). By virtue of this
result, a $\Gamma$-contraction $(S, P)$ is said to be
\textit{c.n.u.} if the contraction $P$ is c.n.u.

Let $P \in \clb(\clh)$ be a contraction and $V \in \clb(\clk)$ be an
isometry. If $V$ is an isometric dilation of $P$, then there exists
an isometry $\Pi : \clh \raro \clk$ such that \[\Pi P^* = V^* \Pi.\]
Conversely, if an isometry $\Pi : \clh \raro \clk$ intertwine $P^*$
and $V^*$, then that $V$ is an isometric dilation of $P$. In the
sequel, we shall identify an isometric dilation of a contraction $P$
by either the dilation map $V$ on the dilation space $\clk$ or by
the isometry $\Pi$ intertwining $P^*$ and $V^*$. In either case, we
call it an \textit{isometric dilation} of the contraction $P$.
Moreover, if the isometric dilation is minimal, that is, if
\[\clk = \overline{\mbox{span}} \{V^m (\Pi \clh) : m \in
\mathbb{N}\},\]then we say that $\Pi$ is a \textit{minimal isometric
dilation} of $P$

We need to recall two dilation results. The $\Gamma$-isometric
dilation of $\Gamma$-contractions is due to Agler and Young
\cite{AY5}, while the isometric dilation of c.n.u. contractions is
due to Sz.-Nagy \cite{Nagy} and Sz.-Nagy and Foias \cite{NF}.

\begin{Theorem}\label{gamma-dilation-AY}
\textsf{(Agler and Young)} Let $(S, P)$ be a $\Gamma$-contraction on
a Hilbert space $\clh$. Then there exists a Hilbert space $\clk$
containing $\clh$, a $\Gamma$-co-isometry $(\tilde{S}, \tilde{P})$
on $\clk$ and an orthogonal decomposition $\clk_1 \oplus \clk_2$ of
$\clk$ such that:

(i) $\clk_1$ and $\clk_2$ are joint invariant subspaces of
$\tilde{S}$ and $\tilde{P}$, and $S = \tilde{S}|_{\clh}, P =
\tilde{P}|_{\clh}$;

(ii) $\clk_1$ and $\clk_2$ reduce both $\tilde{S}$ and $\tilde{P}$;
and

(iii) $(\tilde{S}|_{\clk_1}, \tilde{P}|_{\clk_1})$ is a pure
$\Gamma$-isometry and $(\tilde{S}|_{\clk_2}, \tilde{P}|_{\clk_2})$
is a $\Gamma$-unitary.
\end{Theorem}

\begin{Theorem}\label{Nagy-Foias dilation}\textsf{(Sz.-Nagy)} Let $P$ be a c.n.u. contraction on a Hilbert space $\clh$.
Then there exists a Hilbert space $\clk$ containing $\clh$, an
isometry $V$ on $\clk$ and an orthogonal decomposition $\clk =
\clk_1 \oplus \clk_2$ of $\clk$ such that:

(i) $\clk_1$ and $\clk_2$ are invariant subspaces of $V$ and $P =
V|_{\clh}$;

(ii) $\clk_1$ and $\clk_2$ are reducing subspaces of $V$;

(iii) $V|_{\clk_1}$ is a pure isometry and $V|_{\clk_2}$ is a
unitary; and

(iv) the dilation $V$ on $\clk$ is unique (up to unitary
equivalence) when it is assumed to be minimal.
\end{Theorem}

\NI We like to point out the absence of the minimality property of
the $\Gamma$-isometry in Theorem \ref{gamma-dilation-AY}.

We still need to develop few more definitions and notations. Let $P$
be a contraction on a Hilbert space $\clh$. Then the defect
operators of $P$ are defined by \[D_P = (I_{\clh} - P^*
P)^{\frac{1}{2}} \in \clb(\clh) \quad \mbox{and}\quad D_{P^*} =
(I_{\clh} - P P^*)^{\frac{1}{2}} \in \clb(\clh),\]and the defect
spaces by \[\cld_P = \overline{\mbox{ran} D_P} \quad \mbox{and}\quad
\cld_{P^*} = \overline{\mbox{ran} D_{P^*}}.\]The characteristic
function $\Theta_P \in H^{\infty}_{\clb(\cld_P, \cld_{P^*})}
(\mathbb{D})$ is defined by
\[\Theta_P (z) = [ - P + z D_{P^*} (I_{\clh} - z P^*)^{-1}
D_P]|_{\cld_P}, \quad (z \in \mathbb{D})\]which yields the
multiplication operator $M_{\Theta_P} \in
\clb(H^2_{\cld_P}(\mathbb{D}), H^2_{\cld_{P^*}}(\mathbb{D}))$
defined by \[(M_{\Theta_P} f)(z) = \Theta_P(z) f(z),\]for all $f \in
H^2_{\cld_P}(\mathbb{D})$ and $z \in \mathbb{D}$. Note that
\[M_{\Theta_P} (M_z \otimes I_{\cld_P}) = (M_z \otimes
I_{\cld_{P^*}}) M_{\Theta_P}.\] Define \[\Delta_P(t) = [ I_{\cld_P}
- \Theta_P(e^{it})^* \Theta_P(e^{it})]^{\frac{1}{2}}, \quad (t \in
[0, 1])\]on $L^2_{\cld_P}(\mathbb{T})$ and \[\clh_P =
H^2_{\cld_{P^*}}(\mathbb{D}) \oplus \overline{\Delta_P
L^2_{\cld_P}(\mathbb{T})},\]and the subspace \[\cls_P =
\{M_{\Theta_P} f \oplus \Delta_P f : f \in
H^2_{\cld_P}(\mathbb{D})\} \subseteq \clh_P.\]Notice that $M_z
\oplus M_{e^{it}}|_{\overline{\Delta_P L^2_{\cld_P}(\mathbb{T})}}$
on $\clh_P$ is an isometry where $M_z$ on
$H^2_{\cld_{P^*}}(\mathbb{D})$ is the pure part and
$M_{e^{it}}|_{\overline{\Delta_P L^2_{\cld_P}(\mathbb{T})}}$ on
$\overline{\Delta_P L^2_{\cld_P}(\mathbb{T})}$ is the unitary part
in the sense of the Wold decomposition of isometries. Moreover,
$\cls_P$ is invariant under $M_z \oplus
M_{e^{it}}|_{\overline{\Delta_P L^2_{\cld_P}(\mathbb{T})}}$. Define
the quotient space
\[\mathcal{Q}_P = \clh_P \ominus \cls_P.\]

Let $\cls$ be a closed subspace of a Hilbert space $\clh$. We shall
denote the orthogonal projection from $\clh$ onto $\cls$ by
$\bm{P}_{\cls}$.

\begin{Theorem}\label{Nagy-Foias-model}\textsf{(Sz.-Nagy and
Foias)} Let $P$ be a c.n.u. contraction on a Hilbert space $\clh$.
Then

(i) $P$ is unitarily equivalent to $\bm{P}_{\clq_P}[M_z \oplus
M_{e^{it}}|_{\overline{\Delta_P
L^2_{\cld_P}(\mathbb{T})}}]|_{\clq_P}$.

(ii) The minimal isometric dilation of $P$ can be identified with
$M_z \oplus M_{e^{it}}|_{\overline{\Delta_P
L^2_{\cld_P}(\mathbb{T})}}$ on $\clh_P$.
\end{Theorem}

In the sequel, by virtue of the unitary $U : H^2_{\cle}(\mathbb{D})
\raro H^2(\mathbb{D}) \otimes \cle$ defined by \[z^m \eta \mapsto
z^m \otimes \eta, \quad (\eta \in \cle, m \in \mathbb{N})\]we shall
often identify the vector valued Hardy space
$H^2_{\cle}(\mathbb{D})$ with $H^2(\mathbb{D}) \otimes \cle$.

\vspace{0.2in}

\newsection{Beurling-Lax-Halmos Representations of $\Gamma$-isometries}

This section will focus on a characterization of joint invariant
subspaces of pure $\Gamma$-isometries.

It is well known that the only invariant of pure unweighted
unilateral shift operators is the multiplicity. That is, $M_z$ on
$H^2_{\cle}(\mathbb{D})$ and $M_z$ on $H^2_{\clf}(\mathbb{D})$ are
unitarily equivalent if and only if $\cle$ and $\clf$ are isomorphic
Hilbert spaces. We begin with a characterization of pure
$\Gamma$-isometries in terms of the symbols associated with them.

\begin{Theorem}\label{gamma-iso}
Let $A \in \clb(\cle)$ and $B \in \clb(\clf)$. Then $(M_{A+A^*z},
M_z)$ on $H^2_{\cle}(\mathbb{D})$ and $(M_{B+B^*z}, M_z)$ on
$H^2_{\clf}(\mathbb{D})$ are unitarily equivalent if and only if $A$
and $B$ are unitarily equivalent.
\end{Theorem}
\NI \textsf{Proof.} Let $U : \cle \raro \clf$ be a unitary operator
such that $UA = BU$. Then the unitary operator
\[\tilde{U} = I_{H^2(\mathbb{D})} \otimes U : H^2(\mathbb{D}) \otimes \cle \raro H^2(\mathbb{D}) \otimes
\clf,\]intertwine the corresponding multiplication operators.
Moreover,
\[\tilde{U} (I_{H^2(\mathbb{D})} \otimes A + M_z \otimes A^*) = (I_{H^2(\mathbb{D})} \otimes B + M_z \otimes
B^*) \tilde{U}.\]This proves the sufficiency part.

\NI Conversely, let $\tilde{U} : H^2(\mathbb{D}) \otimes \cle \raro
H^2(\mathbb{D}) \otimes \clf$ be a unitary operator and
\[\tilde{U} (I_{H^2(\mathbb{D})} \otimes A + M_z \otimes A^*) =
(I_{H^2(\mathbb{D})} \otimes B + M_z \otimes B^*) \tilde{U},\]and
\[\tilde{U} (M_z \otimes I_{\cle}) = (M_z \otimes I_{\clf})
\tilde{U}.\]From the last equality it follows that $\tilde{U} =
I_{H^2(\mathbb{D})} \otimes U$ for some unitary operator $U : \cle
\raro \clf$. Then \[ I_{H^2(\mathbb{D})} \otimes U A U^* + M_z
\otimes U A^* U^* = I_{H^2(\mathbb{D})} \otimes B + M_z \otimes
B^*\]implies that $UAU^* = B$, and this completes the proof. \qed

\vspace{0.2in}

The following corollary is a simple but instructive characterization
of pure $\Gamma$-isometries.

\begin{Corollary}
Let $(S_i, P_i)$ be a pair of pure $\Gamma$-isometries on
$H^2_{\cle_i}(\mathbb{D})$ where $i = 1, 2$. Then $(S_1, P_1)$ and
$(S_2, P_2)$ are unitarily equivalent if and only if \[(S^*_1 - S_1
P_1^*) \cong (S^*_2 - S_2 P_2^*).\]
\end{Corollary}

\NI \textsf{Proof.} Let \[(S, P) = (I_{H^2(\mathbb{D})} \otimes A +
M_z \otimes A^*, M_z \otimes I_{\cle}),\] be a pure
$\Gamma$-isometry on $H^2_{\cle}(\mathbb{D})$. Then
\[\begin{split} S^* - S P^* & = (I_{H^2(\mathbb{D})} \otimes A + M_z \otimes A^*)^* - (I_{H^2(\mathbb{D})}
\otimes A + M_z \otimes A^*) (M_z \otimes I_{\cle})^*\\ & =
I_{H^2(\mathbb{D})} \otimes A^* + M_z^* \otimes A - M^*_z \otimes A
- ((I_{H^2(\mathbb{D})} - \bm{P}_{\mathbb{C}}) \otimes A^*)\\& =
\bm{P}_{\mathbb{C}} \otimes A^*,
\end{split}\]where $\bm{P}_{\mathbb{C}}$ is the orthogonal projection
from $H^2(\mathbb{D})$ onto the space of constant functions in
$H^2(\mathbb{D})$. Consequently, by the previous theorem  $S_1^* -
S_1 P_1^*$ and $S_2^* - S_2 P_2^*$ are unitarily equivalent if and
only if $(S_1, P_1)$ and $(S_2, P_2)$ are unitarily equivalent. This
completes the proof. \qed

A closed subspace $\cls \neq \{0\}$ of $H^2_{\cle_*}(\mathbb{D})$ is
said to be \textit{$(A + A^* M_z, M_z)$-invariant} if $\cls$ is
invariant under both $A + A^* M_z$ and $M_z$.

Let $\cls \neq \{0\}$ be a closed subspace of
$H^2_{\cle_*}(\mathbb{D})$. By virtue of the Beurling-Lax-Halmos
theorem $\cls$ is $M_z$-invariant if and only if there exists a
Hilbert space $\cle$ and an inner function $\Theta \in
H^{\infty}_{\clb(\cle, \cle_*)}(\mathbb{D})$ such that \[\cls =
M_{\Theta} H^2_{\cle}(\mathbb{D}).\] Moreover, the pair $\{\cle,
\Theta\}$ is unique in an appropriate sense (cf. see \cite{NF}).

Let $\cls$ be a non-zero $(A + A^* M_z, M_z)$-invariant subspace of
$H^2_{\cle_*}(\mathbb{D})$. Then in particular, by the
Beurling-Lax-Halmos theorem
\[\cls = M_{\Theta} H^2_{\cle}(\mathbb{D}),\]for some Hilbert space
$\cle$ and inner multiplier $\Theta$. The next theorem will show
that $\cls$ is $(A + A^* M_z, M_z)$-invariant if and only
$M_{\Theta}$ intertwine $A + A^* M_z$ and $B + B^* M_z$ for some
unique $B \in \clb(\cle)$ with $w(B) \leq 1$.

Before we proceed with the formal statement and the proof let us
remark that the classification result answers a question left open
by Agler and Young in \cite{AY5}. The proof is remarkably simple (a
straightforward application of Theorem \ref{gamma-isometry}) and may
be of independent interest. However, the intuitive idea behind this
``guess'' is that, $(A + A^* M_z, M_z)$ turns
$H^2_{\cle_*}(\mathbb{D})$ into a natural Hilbert module over
$\mathbb{C}[z_1, z_2]$ (see Section 7).

\begin{Theorem}\label{BLH}
Let $\cls \neq \{0\}$ be a closed subspace of
$H^2_{\cle_*}(\mathbb{D})$ and $A \in \clb(\cle_*)$ with $w(A) \leq
1$. Then $\cls$ is a $(M_{A+A^*z}, M_z)$-invariant subspace if and
only if
\[(A + A^* M_z) M_{\Theta} = M_{\Theta} (B + B^* M_z),\]for some operator $B \in
\clb(\cle)$ with $w(B) \leq 1$, where $\cls = M_{\Theta}
H^2_{\cle}(\mathbb{D})$ is the Beurling-Lax-Halmos representation of
$\cls$. Moreover, when such an operator $B$ exists it is unique (up
to unitary equivalence).
\end{Theorem}

\NI\textsf{Proof.} Let $\cls \neq \{0\}$ be a $(M_{A+A^*z},
M_z)$-invariant subspace and \[\cls = M_{\Theta}
H^2_{\cle}(\mathbb{D}),\]be the Beurling-Lax-Halmos representation
of $\cls$ where $\Theta \in H^{\infty}_{\clb(\cle,
\cle_*)}(\mathbb{D})$ is an inner multiplier and $\cle$ is an
auxiliary Hilbert space. Also
\[(A + A^* M_z) (M_{\Theta} H^2_{\cle}(\mathbb{D})) \subseteq
M_{\Theta} H^2_{\cle}(\mathbb{D}),\]implies that \[(A + A^* M_z)
M_{\Theta} = M_{\Theta} M_{\Psi},\]for some unique $\Psi \in
H^{\infty}_{\clb(\cle)}(\mathbb{D})$. Therefore,\[M_{\Theta}^* (A +
A^* M_z) M_{\Theta} = M_{\Psi}.\]Multiplying both sides by $M_z^*$
we have
\[M_z^* M_{\Theta}^* (A + A^* M_z) M_{\Theta} = M_z^* M_{\Psi}.\]Then
\[M_{\Theta}^* (A M_z^* + A^*) M_{\Theta} = M^*_z
M_{\Psi}.\]Consequently, $M_z^* M_{\Psi} = M_{\Psi}^*$, or
equivalently, $M_{\Psi} = M_{\Psi}^* M_z$. Since $\|M_{\Psi}\| \leq
2$, that $(M_{\Psi}, M_z)$ is a $\Gamma$-isometry, and hence by
Theorem \ref{gamma-isometry}, it follows that
\[M_{\Psi} = B + B^* M_z,\]for some $B \in \clb(\cle)$ and $w(B) \leq
1$, and uniqueness of $B$ follows from that of $\Psi$.

\NI The converse part is trivial, and the proof is complete. \qed

To complete this section we will present the following variant of
our Beurling-Lax-Halmos theorem for $\Gamma$-isometries.

\begin{Theorem}\label{BLH-mult}
Let $\Theta \in H^\infty_{\clb(\cle, \cle_*)}(\mathbb{D})$ be an
inner function and $A \in \clb(\cle_*)$. Then $\cls = M_{\Theta}
H^2_{\cle}(\mathbb{D}) \subseteq H^2_{\cle_*}(\mathbb{D})$ is
invariant under the Toeplitz operator with analytic polynomial
symbol $A + A^*z$ if and only if there exists an operator $B \in
\clb(\cle)$ such that
\[(A + A^* z) \Theta(z) = \Theta(z) (B + B^* z). \quad \quad (z \in \mathbb{D})\]Moreover, when such an operator $B$ exists it is unique (up
to unitary equivalence).
\end{Theorem}

We like to point out that the above result is an application of the
theory of $\Gamma$-contractions to the classical one-variable
operator theory. Moreover, our result suggests a tentative
connection between the theory of spectral sets and invariant
subspaces of Toeplitz operator with analytic polynomial symbol. We
will discuss some of these extensions at the end of this paper.

\vspace{0.2in}

\newsection{Representation of $\Gamma$-contractions}

In this section we will show that a c.n.u. $\Gamma$-contraction can
be realized as a compression of a $\Gamma$-isometry in the Sz.-Nagy
and Foias model space $\clq_P$ of the c.n.u. contraction $P$.
Moreover, we show that the representation of $S$ in $\clq_P$ is
given by a direct sum of a $\Gamma$-isometry and a $\Gamma$-unitary.
Our method involves a "pull-back" technique of the Agler-Young's
isometric dilation to the Sz.-Nagy and Foias minimal isometric
dilation.

First note that if $(S,P)$ is a $\Gamma$-contraction, then $(S^*,
P^*)$ is also a $\Gamma$-contraction, which is equivalent to saying
that $(S^*, P^*)$ has a $\Gamma$-isometric dilation. More precisely,
let $(S,P)$ be a $\Gamma$-contraction on a Hilbert space $\clh$.
Then there exists Hilbert spaces $\clk_1$ and $\clk_2$ and a pure
$\Gamma$-isometry $(\tilde{S}_i, \tilde{P}_i)$ on $\clk_1$ and a
$\Gamma$-unitary $(\tilde{S}_u, \tilde{P}_u)$ on $\clk_2$ and a
isometry (see Theorem \ref{gamma-dilation-AY} or Theorem 3.2 in
\cite{AY5}) \[\Pi_{AY} : \clh \raro \clk_1 \oplus \clk_2,\] such
that\[\Pi_{AY} S^* = (\tilde{S}_i \oplus \tilde{S}_u)^* \Pi_{AY}
\quad \mbox{and}\quad \Pi_{AY} P^* = (\tilde{P}_i \oplus
\tilde{P}_u)^* \Pi_{AY}.\]

\NI On the other hand, by Theorem \ref{Nagy-Foias-model}, we have
the Sz.-Nagy and Foias isometric dilation
\[\Pi_{NF} : \clh \raro \clh_P,\]of the c.n.u. contraction $P$ on $\clh$ with
\[\Pi_{NF} P^* = (M_z \oplus M_{e^{it}}|_{\overline{\Delta_P
L^2_{\cld_P}(\mathbb{T})}})^* \Pi_{NF}.\] Moreover, this dilation is
minimal and hence unique.

The following factorization theorem provides a connection between
the minimal isometric dilation to any other isometric dilation of a
given contraction.

\begin{Theorem}\label{factor}\textsf{(Factorization of Dilations)}
Let $P$ be a c.n.u. contraction on a Hilbert space $\clh$. Let $V$
on $\clk$ be an isometric dilation of $P$ corresponding to an
isometry $\Pi : \clh \raro \clk$. Then there exists a unique
isometry $\Phi \in \clb(\clh_P, \clk)$ such that
\[\Pi = \Phi \Pi_{NF},\]and
\[\Phi (M_z \oplus M_{e^{it}}|_{\overline{\Delta_P L^2_{\cld_P}(\mathbb{T})}})^* = V^* \Phi.\]Moreover,
let $\clk = H^2_{\cle}(\mathbb{D}) \oplus \clk_u$ and $V = M_z
\oplus U$ be the Wold decomposition of $V$ for some unitary $U \in
\clb(\clk_u)$. Then
\[\Phi = (I_{H^2(\mathbb{D})} \otimes V_1) \oplus V_2,\]for some
isometries $V_1 \in \clb(\cld_{P^*}, \cle)$ and $V_2 \in
\clb(\overline{\Delta_P L^2_{\cld_P}(\mathbb{T})}, \clk_u)$.
\end{Theorem}

\NI\textsf{Proof.} Since $\Pi_{NF} : \clh \raro \clh_P$ is the
minimal isometric dilation of $P$ we have
\[\clh_P = \mathop{\bigvee}_{m=0}^{\infty} (M_z \oplus M_{e^{it}}|_{\overline{\Delta_P L^2_{\cld_P}(\mathbb{T})}})^m (\Pi_{NF}
\clh).\]Furthermore, the $V$-reducing subspace
\[\clk_m := \mathop{\bigvee}_{m=0}^{\infty}
V^m (\Pi \clh) \subseteq \clk,\]is the minimal isometric dilation
space of $P$ and hence there exists an isometry \[\Phi : \clh_p
\raro \clk_m \oplus \clk_r,\]defined by
\begin{equation}\label{unique-Phi}\Phi (M_z \oplus
M_{e^{it}}|_{\overline{\Delta_P L^2_{\cld_P}(\mathbb{T})}})^m
(\Pi_{NF} h) = V^m (\Pi h),\end{equation}for all $h \in \clh$ and $m
\in \mathbb{N}$, where $\clk_r = \clk \ominus \clk_m$. Since
\[\begin{split} &\Phi (M_z \oplus
M_{e^{it}}|_{\overline{\Delta_P L^2_{\cld_P}(\mathbb{T})}})^* (M_z
\oplus M_{e^{it}}|_{\overline{\Delta_P L^2_{\cld_P}(\mathbb{T})}})^m
\Pi_{NF}  \\ & = \Phi (M_z \oplus M_{e^{it}}|_{\overline{\Delta_P
L^2_{\cld_P}(\mathbb{T})}})^{m-1} \Pi_{NF} = V^{m-1} \Pi = V^*(V^{m}
\Pi)\\ & = V^* \Phi (M_z \oplus M_{e^{it}}|_{\overline{\Delta_P
L^2_{\cld_P}(\mathbb{T})}})^m \Pi_{NF},
\end{split}\]
for all $m \geq 1$ and \[\Phi (M_z \oplus
M_{e^{it}}|_{\overline{\Delta_P L^2_{\cld_P}(\mathbb{T})}})^*
\Pi_{NF} = \Phi \Pi_{NF} P^* = \Pi P^* = V^* \Pi,\] it follows that
\[\Phi (M_z \oplus M_{e^{it}}|_{\overline{\Delta_P L^2_{\cld_P}(\mathbb{T})}})^* =
V^* \Phi.\]To prove the last part, let $\clk =
H^2_{\cle}(\mathbb{D}) \oplus \clk_u$ and $V = M_z \oplus U$ for
some unitary $U \in \clb(\clk_u)$. Let \[\Phi = \begin{bmatrix}
X_1&X_2\\X_3&X_4\end{bmatrix} : \clh_P =
H^2_{\cld_{P^*}}(\mathbb{D}) \oplus \overline{\Delta_P
L^2_{\cld_P}(\mathbb{T})} \raro \clk = H^2_{\cle}(\mathbb{D}) \oplus
\clk_u.\]Then by the intertwining property of $\Phi$ with the
conjugates of the multiplication operators we have that \[X_1 M_z^*
= M_z^* X_1, \quad X_4 M_{e^{it}}^*|_{\overline{\Delta_P
L^2_{\cld_P}(\mathbb{T})}} = U^* X_4,\] and \[X_2
M_{e^{it}}^*|_{\overline{\Delta_P L^2_{\cld_P}(\mathbb{T})}} = M_z^*
X_2, \quad X_3 M_z^* = U^* X_3.\]Since both $X_2^*$ and $X_3^*$
intertwine a unitary and a pure isometry, it follows that (cf. Lemma
2.5 in \cite{AY5})
\[X_2 = 0 \quad \mbox{and} \quad X_3 = 0.\] Therefore \[\Phi
= \begin{bmatrix} X_1&0\\0&V_2\end{bmatrix},\]where $X_4 = V_2$.
Finally, since \[\mbox{ran} \Phi = \mbox{ran} X_1 \oplus
\mbox{ran}V_2 \subseteq H^2_{\cle}(\mathbb{D}) \oplus \clk_u,\] is a
$(M_z \oplus U)$-reducing subspace of $\clk$ and $\mbox{ran} X_1
\subseteq H^2_{\cle}(\mathbb{D})$, it follows that $\mbox{ran} X_1$
is a $M_z$-reducing subspace of $H^2_{\cle}(\mathbb{D})$.
Consequently,
\[X_1 = I_{H^2(\mathbb{D})} \otimes V_1,\]for some isometry
$V_1 \in \clb(\cld_{P^*}, \cle)$.

\NI Uniqueness of $\Phi$ follows from the equality
(\ref{unique-Phi}). This completes the proof. \qed

\vspace{0.2in}

The above factorization result can be summarized in the following
commutative diagram:

 \setlength{\unitlength}{3mm}
 \begin{center}
 \begin{picture}(40,16)(0,0)
\put(15,3){$\clh$}\put(19,1.6){$\Pi$} \put(22.9,3){$\clk$} \put(22,
10){$\clh_P$} \put(22.8,9.2){ \vector(0,-1){5}} \put(15.8,
4.3){\vector(1,1){5.8}} \put(16.4,
3.4){\vector(1,0){6}}\put(16.5,8){$\Pi_{NF}$}\put(24.3,7){$\Phi$}
\end{picture}
\end{center}
where $\Phi$ is a unique isometry which intertwines the adjoints of
the multiplication operators. Related results along this line can be
found in the context of the commutant lifting theorem of
contractions (cf. page 134 in \cite{FF} and page 133 in
\cite{FFGK}). However, the minimal isometric dilation space in this
consideration is the Schaffer's dilation space.

Before we proceed further let us recall the \textit{Putnam's
Corollary} \cite{Pu}: If $A$ and $B$ are normal operators on Hilbert
spaces $\clh$ and $\clk$ respectively and if $C$ is a bounded linear
operator from $\clh$ to $\clk$ such that $CA = BC$, then $CA^* = B^*
C$.

The main result of this section is the following theorem concerning
an analytic model of a c.n.u. $\Gamma$-contraction.

\begin{Theorem}\label{gamma-dilation}
Let $(S, P)$ be a c.n.u. $\Gamma$-contraction on a Hilbert space
$\clh$. Then

\begin{equation}\label{*}\tag{\textsf{NF-AY}}S \cong \bm{P}_{\clq_P} ((A+A^*M_z) \oplus U)|_{\clq_P},\end{equation}
where $A \in \clb(\cld_{P^*})$ with $w(A) \leq 1$ and $U$ is in
$\clb(\overline{\Delta_P L^2_{\cld_P}(\mathbb{T})})$ such that \[(U,
M_{e^{it}}|_{\overline{\Delta_P L^2_{\cld_P}(\mathbb{T})}}),\]is a
$\Gamma$-unitary.
\end{Theorem}

\NI \textsf{Proof.} Let \[\Pi_{NF} : \clh \raro \clh_P =
H^2_{\cld_{P^*}}(\mathbb{D}) \oplus \overline{\Delta_P
L^2_{\cld_P}(\mathbb{T})},\] be the Sz.-Nagy and Foias minimal
isometric dilation of the c.n.u. contraction $P$ as in Theorem
\ref{Nagy-Foias-model} and
\[\Pi_{AY} : \clh \raro \clk_1 \oplus \clk_2,\]be the Agler-Young's
$\Gamma$-isometric dilation of the c.n.u. $\Gamma$-contraction $(S,
P)$ as in Theorem \ref{gamma-dilation-AY}, where $\clk_1 =
H^2_{\cle}(\mathbb{D})$ is the pure part and $\clk_2$ is the unitary
part. By Theorem \ref{factor}, we have the following commutative
diagram

\setlength{\unitlength}{3mm}
\begin{center}
\begin{picture}(40,16)(0,0)
\put(15,3){$\clh$}\put(19,1.6){$\Pi_{AY}$} \put(22.3,3){$\clk_1
\oplus \clk_2$} \put(23, 10){$\clh_P$} \put(23.8,9.2){
\vector(0,-1){5}} \put(15.8, 4.3){\vector(1,1){5.8}} \put(16.2,
3.4){\vector(1,0){5.8}}\put(16.5,8){$\Pi_{NF}$}\put(25.3,7){$\tilde{V}$}
\end{picture}
\end{center}

\NI where $\tilde{V}$ is an isometry of the form \[\tilde{V}  =
(I_{H^2(\mathbb{D})} \otimes \hat{V_1}) \oplus \hat{V_2},\]for some
isometries $\hat{V_1} \in \clb(\cld_{P^*}, \cle)$ and $\hat{V_2} \in
\clb(\overline{\Delta_P L^2_{\cld_P}(\mathbb{T})}, \clk_2)$.
Moreover,
\begin{equation}\label{V2-int}\hat{V_2} M_{e^{it}}^*|_{\overline{\Delta_P
L^2_{\cld_P}(\mathbb{T})}} = \tilde{P}_u^* \hat{V_2}.\end{equation}
Then
\[\Pi_{AY} = ((I_{H^2(\mathbb{D})} \otimes \hat{V_1}) \oplus \hat{V_2}) \Pi_{NF}.\]
Since \[ \Pi_{AY} P^* = (\tilde{P}_i \oplus \tilde{P}_u)^* \Pi_{AY}
= ((M_z \otimes I_{\cle}) \oplus \tilde{P}_u)^* \Pi_{AY},\] and
\[\Pi_{AY} S^* = (\tilde{S}_i \oplus \tilde{S}_u)^* \Pi_{AY} = ((I_{H^2(\mathbb{D})}
\otimes A + M_z \otimes A^*) \oplus \tilde{S}_u)^* \Pi_{AY},\] we
have
\begin{equation}\label{P*}
P^* = \Pi_{NF}^* ((M_z \otimes I_{\cle}) \oplus  \hat{V_2}^*
\tilde{P_u} \hat{V_2})^* \Pi_{NF}, \end{equation} and

\begin{equation}\label{S*}
S^* = \Pi_{NF}^* ((I_{\cle} \otimes \hat{V_1}^* A \hat{V_1} + M_z
\otimes \hat{V_1}^* A^*\hat{V_1}) \oplus \hat{V_2}^* \tilde{S}_u
\hat{V_2})^* \Pi_{NF}.
\end{equation}
By (\ref{V2-int}) and Putnam's Corollary, we have
\[\hat{V}_2  M_{e^{it}}|_{\overline{\Delta_P
L^2_{\cld_P}(\mathbb{T})}} = \tilde{P}_u \hat{V}_2.\]In particular,
$\mbox{ran} \hat{V}_2$ is a $\tilde{P}_u$-reducing subspace, and
\[\hat{V}_2^* \tilde{P}_u \hat{V}_2 =
M_{e^{it}}|_{\overline{\Delta_P L^2_{\cld_P}(\mathbb{T})}} \in
\clb(\overline{\Delta_P L^2_{\cld_P}(\mathbb{T})}).\] Consequently,
\[\begin{split}\Pi_{NF} P^* \Pi_{NF}^* & = \bm{P}_{\clq_P} ((M_z \otimes I_{\cld_{P^*}})
\oplus  \hat{V_2}^* \tilde{P_u} \hat{V_2})^*|_{\clq_P} \\ & =
\bm{P}_{\clq_P} ((M_z \otimes I_{\cld_{P^*}}) \oplus
M_{e^{it}}|_{\overline{\Delta_P
L^2_{\cld_P}(\mathbb{T})}})^*|_{\clq_P},\end{split}\] and \[\Pi_{NF}
S^* \Pi_{NF}^* = \bm{P}_{\clq_P} ((I_{\cld_{P^*}} \otimes
\hat{V_1}^* A \hat{V_1} + M_z \otimes \hat{V_1}^* A^*\hat{V_1})
\oplus \hat{V_2}^* \tilde{S}_u \hat{V_2})^*|_{\clq_P}.
\]Therefore, \[S^* \cong \bm{P}_{\clq_P} ((\tilde{A} +\tilde{A}^*M_z) \oplus
\tilde{U})|_{\clq_P},\]where $\tilde{A} = \hat{V_1}^* A \hat{V_1}$
and $\tilde{U} =  \hat{V_2}^* \tilde{S_u} \hat{V_2}$. Since $w(A)
\leq 1$ we have that $w(\tilde{A}) \leq 1$.

\NI It remains to prove that \[(\tilde{U}, \hat{V_2}^* \tilde{P_u}
\hat{V_2}) = (\tilde{U},  M_{e^{it}}|_{\overline{\Delta_P
L^2_{\cld_P}(\mathbb{T})}})\]is a $\Gamma$-unitary. Since
$(\tilde{S_u}, \tilde{P_u})$ is a $\Gamma$-unitary, we conclude that
\[\tilde{S_u} = \tilde{S_u}^*\tilde{P_u}.\] Using the fact that the
range of $\hat{V_2}$ is $\tilde{P_u}$-reducing, it follows that
\[\tilde{S_u} \hat{V_2} = \tilde{S_u}^* \tilde{P_u} \hat{V_2} =
(\tilde{S_u}^* \hat{V_2}) (\hat{V_2}^* \tilde{P_u}
\hat{V_2}).\]Hence
\[\hat{V_2}^*\tilde{S_u} \hat{V_2} = (\hat{V_2}^* \tilde{S_u}
\hat{V_2})^* (\hat{V_2}^* \tilde{P_u} \hat{V_2}),\] which implies
that
\[(\hat{V_2}^* \tilde{S_u} \hat{V_2}, \hat{V_2}^* \tilde{P_u}
\hat{V_2}) = (\tilde{U}, \hat{V_2}^* \tilde{P_u} \hat{V_2}) =
(\tilde{U}, M_{e^{it}}|_{\overline{\Delta_P
L^2_{\cld_P}(\mathbb{T})}})\]is a $\Gamma$-unitary. This finishes
the proof of the theorem. \qed

\vspace{0.1in}

The following result shows that the representation of $S - P S^*$ in
the Sz.-Nagy and Foias quotient space $\clq_P$ is the compression of
the scalar operator $A$ on $\clq_P$.

\begin{Corollary}\label{equal}
With notations as in Theorem \ref{gamma-dilation}, representation of
the operator $S^* - S P^*$ in $\clq_P$ is given by  \[S^* - S P^*
\cong \bm{P}_{\clq_P}(\bm{P}_{\mathbb{C}} \otimes A^*)|_{\clq_P}.\]
\end{Corollary}

\NI \textsf{Proof.} Since
\[\begin{split} &[\bm{P}_{\clq_P} ((A + A^* M_z) \oplus U)|_{\clq_P}]^* -
[\bm{P}_{\clq_P} ((A + A^* M_z) \oplus U)|_{\clq_P}]
\bm{P}_{\clq_P}([M_z \oplus M_{e^{it}})|_{\clq_P}]^*\\ =  &
\bm{P}_{\clq_P}[ ((A^* + A M^*_z) \oplus U^*)
- ((A + A^* M_z) \oplus U) (M_z^* \oplus M^*_{e^{it}})]_{\clq_P} \\
= & \bm{P}_{\clq_P}(\bm{P}_{\mathbb{C}} \otimes A^*)|_{\clq_P},
\end{split}\] we have \[S^* - S P^*
\cong \bm{P}_{\clq_P}(\bm{P}_{\mathbb{C}} \otimes A^*)|_{\clq_P}.\]
This completes the proof. \qed

\vspace{0.3in}

\newsection{Unique representation of $\Gamma$-contractions}

In this section we will discuss the uniqueness of the
$\Gamma$-isometry and the $\Gamma$-unitary part in the
representation (\ref{*}) of a c.n.u. $\Gamma$-contraction.

We begin by recalling one way to construct the minimal isometric
dilation of a c.n.u. contraction. More details can be found, for
instance, in the monograph by Foias and Frazho (Page 137 in
\cite{FF}). Let $P \in \clb(\clh)$ be a c.n.u. contraction. Then
\[\tilde{X}_P := SOT-\lim_{m \raro \infty} P^m P^{*m},\]is a positive
operator on $\clh$. Let $X_P$ be the positive square root of
$\tilde{X}_P$. Then \[\|X_P h\|^2 = \lim_{m\raro
\infty}\|P^{*m}h\|^2,\] and \[\|X_P h\| = \|X_P P^*h\|,\]for all $h
\in \clh$. Consequently, there exists an isometry $V_1 \in
\clb(\overline{X_P \clh})$ such that \[V_1 X_P = X_P P^*.\]Let $V_2$
on $\clk_u$ be the minimal unitary extension of $V_1$ so that \[V_2
X_P = X_P P^*.\] Define $\Pi : \clh \raro
H^2_{\cld_{P^*}}(\mathbb{D}) \oplus \clk_u$ by \[\Pi h = D_{P^*}(I -
zP^*)^{-1} h \oplus X_P h,\]for all $h \in \clh$. Then $\Pi$ is an
isometry and \[\Pi P^* = (M_z^* \oplus U^*) \Pi,\]where $U = V_2^*$.
Moreover, $\Pi$ is minimal and
\[\Pi^* ((\mathbb{S}_w \otimes \eta) \oplus 0) = (I - \bar{w} P)^{-1}
D_{P^*} \eta,\]for all $\eta \in \cld_{P^*}$, where $\mathbb{S}$ is
the Szeg\"{o} kernel on the open unit disk defined by \[\mathbb{S}_w
(z) = (1 - z \bar{w})^{-1},\]for all $z, w \in \mathbb{D}$.

In the proof of the following theorem, we shall identify (by virtue
of Theorem \ref{factor} where $\Phi$ is a unitary) the minimal
isometric dilation $\Pi_{NF}$ with the one described above.

\begin{Theorem}\label{Aequiv}
Let $P$ be a c.n.u. contraction on a Hilbert space $\clh$ and $A \in
\clb(\cld_{P^*})$. Then \[D_{P^*} A D_{P^*} \cong \bm{P}_{\clq_P}
(\bm{P}_{\mathbb{C}} \otimes A)|_{\clq_P}.\]
\end{Theorem}

\NI\textsf{Proof.} Let $\bm{ev}_0 : (H^2(\mathbb{D}) \otimes
\cld_{P^*}) \oplus \overline{\Delta_P L^2_{\cld_P}(\mathbb{T})}
\raro \cld_{P^*}$ be the evaluation operator defined by
\[\bm{ev}_0 (f \oplus g) = f(0),\]for all $f \oplus g \in (H^2(\mathbb{D}) \otimes
\cld_{P^*}) \oplus \overline{\Delta_P L^2_{\cld_P}(\mathbb{T})}$.
Then \[\bm{ev}_0 \Pi_{NF} h = (D_{P^*}(I - z P^*)^{-1} h)(0) =
D_{P^*}h,\]for all $h \in \clh$. From this we readily obtain
\[\bm{ev}_0 \Pi_{NF} = D_{P^*}.\]Moreover, \[\Pi_{NF}^*
(\bm{P}_{\mathbb{C}} \otimes A) \eta = \Pi_{NF}^*(1 \otimes A \eta)
= \Pi_{NF}^*((\mathbb{S}_0 \otimes A \eta) \oplus 0) = D_{P^*} A
\eta,\]for all $\eta \in \cld_{P^*}$. Thus,
\[\Pi_{NF}^*(\bm{P}_{\mathbb{C}} \otimes A) = D_{P^*} A \bm{ev}_0.\]
Consequently, \[\Pi_{NF}^* (\bm{P}_{\mathbb{C}} \otimes A) \Pi_{NF}
= D_{P^*} A \bm{ev}_0 \Pi_{NF} = D_{P^*} A D_{P^*}.\]Then the result
follows from the fact that
\[\Pi_{NF} \Pi_{NF}^* (\bm{P}_{\mathbb{C}} \otimes A) \Pi_{NF} \Pi_{NF}^*
= \Pi_{NF} (D_{P^*} A D_{P^*}) \Pi_{NF}^*.\] \qed

The following corollary is immediate.

\begin{Corollary}\label{5.2}
Let $P$ be a c.n.u. contraction and $A \in \clb(\cld_{P^*})$. Then
$A = 0$ if and only if \[\bm{P}_{\clq_P} (\bm{P}_{\mathbb{C}}\otimes
A)|_{\clq_P} = 0.\]
\end{Corollary}

Now we are ready to prove the uniqueness of the $\Gamma$-isometric
part in (\ref{*}).

\begin{Theorem}\label{unique-A}
Let $(S, P)$ be a c.n.u. $\Gamma$-contraction on a Hilbert space
$\clh$. Then the operator $A \in \clb(\cld_{P^*})$ in the
representation (\ref{*}) is unique. That is, if
\[S \cong \bm{P}_{\clq_P} ((\tilde{A} + \tilde{A}^* M_z) \oplus
\tilde{U})|_{\clq_P},\]where $(\tilde{U},
M_{e^{it}}|_{\overline{\Delta_P L^2_{\cld_P}(\mathbb{T})}})$ is a
$\Gamma$-unitary and $\tilde{A} \in \clb(\cld_{P^*})$, $w(\tilde{A})
\leq 1$ and $\tilde{U} \in \clb(\overline{\Delta_P
L^2_{\cld_P}(\mathbb{T})})$, then $A = \tilde{A}$.
\end{Theorem}

\NI\textsf{Proof.} Let \[\bm{P}_{\clq_P} ((A + A^* M_z) \oplus
U)|_{\clq_P} = \bm{P}_{\clq_P} ((\tilde{A} + \tilde{A}^* M_z) \oplus
\tilde{U})|_{\clq_P}\] where $(\tilde{U},
M_{e^{it}}|_{\overline{\Delta_P L^2_{\cld_P}(\mathbb{T})}})$ is a
$\Gamma$-unitary and $\tilde{A} \in \clb(\cld_{P^*})$, $w(\tilde{A})
\leq 1$ and $\tilde{U} \in \clb(\overline{\Delta_P
L^2_{\cld_P}(\mathbb{T})})$. By Corollary \ref{equal}, we have
\[\bm{P}_{\clq_P}(\bm{P}_{\mathbb{C}} \otimes A)|_{\clq_P} =
\bm{P}_{\clq_P}(\bm{P}_{\mathbb{C}} \otimes
\tilde{A})|_{\clq_P}.\]This and Corollary \ref{5.2} implies that \[A
= \tilde{A}.\] This completes the proof. \qed

The following result plays an important role in the proof of the
uniqueness of $\tilde{U}$ in (\ref{*}).

\begin{Proposition}\label{X=0}
Let $X$ be a bounded linear operator on $\overline{\Delta_P
L^2_{\cld_P}(\mathbb{T})}$ where $X M_{e^{it}}|_{\overline{\Delta_P
L^2_{\cld_P}(\mathbb{T})}} = M_{e^{it}}|_{\overline{\Delta_P
L^2_{\cld_P}(\mathbb{T})}} X$ and $X|_{\clq_P} = 0$. Then $X = 0$.
\end{Proposition}

\NI \textsf{Proof.} Let $X|_{\clq_P} = 0$. Then
\[\mbox{ran} X^* = \mbox{ran} (0 \oplus X^*)  \subseteq \clq_P^{\perp} =
\overline{\mbox{span}} \{\Theta_P f \oplus \Delta_P f : f \in
H^2_{\cld_P}(\mathbb{D})\}.\] If $\Delta_P f \in \mbox{ran} X^*$ for
some $f \in H^2_{\cld_P}(\mathbb{D})$ then $\Theta_P f = 0$, or
equivalently, $\Theta_P^* \Theta_P f = 0$. Therefore,
\[\Delta_P^2 f = f,\]and hence \[\Delta_P (\mbox{ran}
X^*) \subseteq H^2_{\cld_P}(\mathbb{D}).\] Also by \[X
M_{e^{it}}|_{\overline{\Delta_P L^2_{\cld_P}(\mathbb{T})}} =
M_{e^{it}}|_{\overline{\Delta_P L^2_{\cld_P}(\mathbb{T})}} X,\] we
conclude that $\overline{\Delta_P (\mbox{ran} X^*)} \subseteq
H^2_{\cld_P}(\mathbb{D})$ is a $M_{e^{it}}|_{\overline{\Delta_P
L^2_{\cld_P}(\mathbb{T})}}$-reducing subspace. Consequently,
\[\Delta_P (\mbox{ran} X^*) = \{0\}.\] Since $X \in
\clb(\overline{\Delta_P L^2_{\cld_P}(\mathbb{T})})$, we have $X =
0$, which completes the proof. \qed

From the previous proposition we readily obtain the desired
uniqueness of $\tilde{U}$ in (\ref{*}).

\begin{Corollary}
Let $(S, P)$ be a c.n.u. $\Gamma$-contraction on a Hilbert space
$\clh$. Then the operator $U \in \clb(\overline{\Delta_P
L^2_{\cld_P}(\mathbb{T})})$ in the representation of $S$ in
(\ref{*}) is unique.
\end{Corollary}
\NI \textsf{Proof.}  Let \[\bm{P}_{\clq_P} ((A + A^* M_z) \oplus
U)|_{\clq_P} = \bm{P}_{\clq_P} ((A + A^* M_z) \oplus
\tilde{U})|_{\clq_P},\]  where $(\tilde{U},
M_{e^{it}}|_{\overline{\Delta_P L^2_{\cld_P}(\mathbb{T})}})$ is a
$\Gamma$-unitary and $\tilde{A} \in \clb(\cld_{P^*})$, $w(\tilde{A})
\leq 1$ and $\tilde{U} \in \clb(\overline{\Delta_P
L^2_{\cld_P}(\mathbb{T})})$. Then
\[\bm{P}_{\clq_P} [0 \oplus (U - \tilde{U})]|_{\clq_P} = 0.\] By
Proposition \ref{X=0} with $X = (U - \tilde{U})^*$, we have \[U =
\tilde{U}.\] This completes the proof. \qed

Combining the above corollary with Theorems \ref{gamma-dilation} and
\ref{unique-A}, we obtain the unique representation of a c.n.u.
$\Gamma$-contraction $(S, P)$ in the model space $\clq_P$.

\begin{Theorem}\label{lift}
Let $(S, P)$ be a c.n.u. $\Gamma$-contraction on a Hilbert space
$\clh$. Then the representing operators $A$ and $U$ in (\ref{*}) are
unique. That is, if \[S \cong \bm{P}_{\clq_P} ((\tilde{A} +
\tilde{A}^* M_z) \oplus \tilde{U})|_{\clq_P},\] where $(\tilde{U},
M_{e^{it}}|_{\overline{\Delta_P L^2_{\cld_P}(\mathbb{T})}})$ is a
$\Gamma$-unitary and $\tilde{A} \in \clb(\cld_{P^*})$, $w(\tilde{A})
\leq 1$ and $\tilde{U} \in \clb(\overline{\Delta_P
L^2_{\cld_P}(\mathbb{T})})$. Then $A \cong  \tilde{A}$ and $U \cong
\tilde{U}$. Moreover, \[((A + A^* M_z) \oplus U, M_z \oplus
M_{e^{it}}|_{\overline{\Delta_P L^2_{\cld_P}(\mathbb{T})}}),\]is the
minimal isometric dilation of the $\Gamma$-contraction $(S, P)$.
\end{Theorem}

Therefore, a c.n.u. $\Gamma$-contraction $(S, P)$ on $\clh$ is
uniquely determined by $\Theta_P$ (and hence by $\clq_P$) and the
\textit{representing operators} $A$ and $U$.

\vspace{0.3in}

\newsection{A Characterization of $\Gamma$-contractions}

Let $(s, p) \in \mathbb{C}^2$. Then $(s, p)$ is in the symmetrized
bidisc $\Gamma$ if and only if $|p| \leq 1$ and that
\[s = \beta + p \bar{\beta},\] for some $\beta \in \mathbb{C}$ such
that $|\beta| \leq 1$ (see \cite{AY1}). In this section we
generalize the scalar characterization of elements in $\Gamma$ to
the class of $\Gamma$-contractions on Hilbert spaces.

We begin by recalling the Schaffer isometric dilation of a
contraction $P$ on $\clh$. In this case, the dilation space is
defined by $\clk_P := \clh \oplus H^2_{\cld_{P}}(\mathbb{D})$. Let
\[V_P = \begin{bmatrix}P& 0\\\bm{D_{P}}&M_z\end{bmatrix},\]where
$\bm{D_{P}} : \clh \raro H^2_{\cld_{P}}(\mathbb{D})$ is the constant
function defined by \[(\bm{D_{P}} h)(z) = D_{P}h,\]for all $h \in
\clh$ and $z \in \mathbb{D}$. That is, \[V_P(h \oplus f) = P h
\oplus (D_{P}h + M_z f),\]for all $h \oplus f \in \clk_P$. Then
$V_P$ is an isometry and the map
\[\Pi_{Sc} : \clh \raro \clk_P = \clh \oplus
H^2_{\cld_{P}}(\mathbb{D}),\] defined by \[\Pi_{Sc} h = h \oplus
0,\] for all $h \in \clh$, satisfies \[\Pi_{Sc} P^* = V^*_P
\Pi_{Sc}.\]The isometric dilation $\Pi_{Sc}$ is known as the
\textit{Schaffer dilation} of the contraction $P$.

The following result summarizes Theorems 4.2, 4.3 and 4.4 in
\cite{BPR}. For completeness and the reader's convenience, we supply
a proof. Moreover, our view is slightly different and the proof is
considerably short and simple.

\begin{Theorem}\label{A-equation}
Let $(S, P)$ be a $\Gamma$-contraction. Then the Schaffer dilation
of $P$ satisfies \[\Pi_{Sc} S^* = W_A^* \Pi_{Sc},\]for some
$\Gamma$-isometry $(W_A, V_P)$ which is uniquely determined by the
operator $A \in \clb(\cld_{P})$ such that $S - S^* P = D_{P} A
D_{P}$ and $w(A) \leq 1$. Conversely, let $(S, P)$ be a commuting
tuple where $P$ is a contraction and $\|S\| \leq 2$ and $S - S^* P =
D_{P} A D_{P}$ for some $A \in \clb(\cld_{P})$ with $w(A) \leq 1$.
Then $(S,P)$ is a $\Gamma$-contraction.
\end{Theorem}

\NI \textsf{Proof.} Let $(S, P)$ be a $\Gamma$-contraction. First,
we assume that $(S, P)$ is c.n.u. Therefore $P$ is a c.n.u.
contraction and by the factorization of dilations, Theorem
\ref{factor}, we have an isometry $\Phi : \clh_P \raro \clk_P$ such
that \[\Pi_{Sc} = \Phi \Pi_{NF}.\]The Schaffer dilation $\Pi_{Sc}$
is minimal means that $\Phi$ is unitary. Then \[(W, V) := (\Phi((A +
A^* M_z) \oplus U) \Phi^*, \Phi (M_z \oplus
M_{e^{it}}|_{\overline{\Delta_P
L^2_{\cld_P}(\mathbb{T})}})\Phi^*),\]is a $\Gamma$-isometry on
$\clk_P$ and
\[\Pi_{Sc} P^* = V^* \Pi_{Sc}, \quad \quad \mbox{and}\quad \quad
\Pi_{Sc} S^* = W^* \Pi_{Sc}.\]Finally, by taking the orthogonal
direct sum of the unitary part with the c.n.u. part, it follows that
$\Pi_{Sc}$ is the minimal isometric dilation of the
$\Gamma$-contraction $(S, P)$. By the uniqueness of the Schaffer
dilation of contractions we therefore identify that $V$ with $V_P$.

Next we will show that \begin{equation}\label{=W}W =
\begin{bmatrix}S & 0\\ \bm{A^* D_{P}}& A + A^*
M_z\end{bmatrix},\end{equation}for some $A \in \clb(\cld_{P^*})$
with $w(A) \leq 1$. To see this, assume \[W = \begin{bmatrix} S &
0\\W_3& W_4\end{bmatrix},\] and compute \[W^* V_P =
\begin{bmatrix}S^* P + W^*_3 D_{P} & W^*_3 M_z\\ W^*_4 \bm{D_{P}}&
W^*_4 M_z\end{bmatrix}.\]Since $(W, V_P)$ is a $\Gamma$-isometry, we
have $W^* V_P = W$ and so
\begin{equation}\label{W=}\begin{bmatrix}S^* P + W^*_3 D_{P} & W^*_3
M_z\\ W^*_4 \bm{D_{P}}& W^*_4 M_z\end{bmatrix} =
\begin{bmatrix}S & 0\\ W_3& W_4\end{bmatrix}.\end{equation}By $W^*_4 M_z =
W_4$ and that $\|W_4\| \leq 2$ we have \[W_4 = A+ A^* M_z,\]for some
$A \in \clb(\cld_{P})$ and $w(A) \leq 1$. Also
\[W_3 = W^*_4 \bm{D_{P}} = (A + A^* M_z)^* \bm{D_{P}} = \bm{A^*
D_{P}},\]which yields the desired representation of $W$. In the
above equality we used the fact that \[(A M_z^* \bm{D_{P}}h)(z)  =
M^*_z D_{P}h = 0,\]and \[(A^* \bm{D_{P}}h)(z) = A^* D_{P}h,\]for all
$h \in \clh$ and $z \in \mathbb{D}$.

Now we will show that $A$ is uniquely determined by $(S, P)$. For
that, equating the $(1,1)$-th entries in (\ref{W=}) we have \[S^* P
+ W^*_3 D_{P} = S.\]Hence
\[S - S^* P = W^*_3 D_{P^*} = D_{P} A D_{P},\]and that $A$ is
uniquely determined by $(S, P)$.

Therefore, if $(S, P)$ is a $\Gamma$-contraction on $\clh$ then
there exists a unique $A \in \clb(\cld_{P^*})$ with $w(A) \leq 1$
and $S - S^* P = D_{P} A D_{P}$ such that $\Pi_{Sc} : \clh \raro
\clk$ satisfies
\[\Pi_{Sc} S^* = W_A^* \Pi_{Sc},\]where $W_A$ is the operator matrix
in (\ref{=W}).

On the other hand, given a commuting tuple $(S, P)$ on $\clh$, where
$P$ is a contraction and $\|S\| \leq 2$ and $S - S^* P = D_{P} A
D_{P}$ for some $B \in \clb(\cld_{P^*})$ (and hence, unique) the
Schaffer dilation of $P$ satisfies \[\Pi_{Sc} S^* = W_A^*
\Pi_{Sc},\]where $W_A$ is the operator matrix in (\ref{=W}).
Moreover, by the given conditions, it is easy to check that $W_A^*
V_P = W_A$. Since $\|W_A\| = r(W_A) \leq 2$, (see page 598 in
\cite{BPR}) we obtain that $(W_A, V_P)$ is a $\Gamma$-isometry, that
is, $\Pi_{Sc}$ is a $\Gamma$-isometric dilation of $(S, P)$. Here we
are using the fact that if $W$ on $\clh$ commutes with an isometry
$V$ and $W^* V = W$ then $W$ is a hyponormal operator and that $r(W)
= \|W\|$ (see Theorem 1 in \cite{St}). Consequently, $(S, P)$ is a
$\Gamma$-contraction. This completes the proof. \qed

\vspace{0.1in}

Let us remark that the equality $S - S^* P = D_P A D_P$ of a
$\Gamma$-contraction $(S, P)$ also follows by applying Corollary
\ref{equal} and Theorem \ref{Aequiv} and the Wold decomposition
theorem of $\Gamma$-contractions to the $\Gamma$-contraction $(S^*,
P^*)$.

Let $(S, P)$ be a $\Gamma$-unitary. Since the only way to obtain a
$\Gamma$-unitary is to symmetrize a pair of commuting unitary
operators, say $U$ and $U_1$, we let $S = U + U_1$ and $P = U U_1$.
Then $U_1 = U^*P$ and hence $S = U + U^* P$. Therefore, a
$\Gamma$-unitary can be represented by $(U_1 + U_1^* U, U)$ for some
commuting unitary operators $U_1$ and $U$ (see Theorem 2.5 in
\cite{BPR}).

Let $\cle$ be a Hilbert space and $U$ be a unitary operator on some
Hilbert space $\clk$ and $\clq$ be a $(M_z \oplus U)^*$-invariant
subspace of $H^2_{\cle}(\mathbb{D}) \oplus \clk$. A bounded linear
operator $X$ on $\clq$ is said to be \textit{compressed scalar} if
\[X = P_{\clq}(A \oplus U_1)|_{\clq},\]for some $A \in \clb(\cle)$
and a unitary $U_1 \in \clb(\clk)$ such that $U U_1 = U_1 U$. Given
a compressed scalar operator $X = P_{\clq}(A \oplus U_1)|_{\clq}$ on
$\clq$, the \textit{decompressed numerical radius} of $X$, denoted
by $w_{\clq}(X)$, is defined by $w(A)$, the numerical radius of the
symbol operator $A$.

The realization of $\Gamma$-unitary along with the
$\Gamma$-isometric dilation theorem, Theorem \ref{lift}, and Theorem
\ref{A-equation} enables us to prove the following characterization
of $\Gamma$-contractions. Here, and in what follows, if $P$ on
$\clh$ is a c.n.u. contraction then we will assume that $\clh =
\clq_{P}$, the Nagy-Foias quotient space.

\begin{Theorem}\label{XP}
Let $P$ be a c.n.u. contraction and $S \in \clb(\clq_P)$ commutes
with $P$ and $\|S\| \leq 2$. Then $(S, P)$ is a $\Gamma$-contraction
on $\clq_P$ if and only if \[S = X + P X^*,\]for some compressed
scalar operator $X \in \clb(\clq_P)$ with $w_{\clq_P}(X) \leq 1$.
\end{Theorem}

\NI\textsf{Proof.} Let $(S, P)$ be a c.n.u. $\Gamma$-contraction and
let
\[((A + A^* M_z) \oplus (U_1 + U_1^* U), M_z \oplus U)\]on $\clh_P$
be the $\Gamma$-isometric dilation of $(S, P)$ (see Theorem
\ref{lift}), where $U = M_{e^{it}}|_{\overline{\Delta_P
L^2_{\cld_P}(\mathbb{T})}}$ and $w(A) \leq 1$. Now
\[\begin{split}(A + A^* M_z) \oplus (U_1 + U_1^* U) & = (A \oplus
U_1) + (A^* M_z \oplus U_1^* U)\\ & = (A \oplus U_1) + (M_z \oplus
U) (A \oplus U_1)^*.
\end{split}\] Identifying $(S, P)$ with

\[(\bm{P}_{\clq_P}((A + A^* M_z) \oplus ( U_1 + U_1^* U))|_{\clq_P}, \bm{P}_{\clq_P} (M_z \oplus U)|_{\clq_P}),\]
on the quotient space $\clq_P$, we have
\[\begin{split} S^* & = ((A + A^* M_z) \oplus (U_1 + U_1^*
U))^*|_{\clq_P}\\ & = ((A \oplus U_1)^* + (A \oplus U_1) (M_z \oplus
U)^*)|_{\clq_P}\\& = \bm{P}_{\clq_P}(A \oplus U_1)^*|_{\clq_P} +
\bm{P}_{\clq_P} (A \oplus U_1) \bm{P}_{\clq_P} (M_z \oplus
U)^*|_{\clq_P}\\& = \bm{P}_{\clq_P}(A \oplus U_1)^*|_{\clq_P} +
\bm{P}_{\clq_P} (A \oplus U_1)|_{\clq_P} \bm{P}_{\clq_P} (M_z \oplus
U)^*|_{\clq_P}\\& = X^* + X P^*,
\end{split}
\]that is, \[S = X + P X^*,\] where \[X := \bm{P}_{\clq_P} (A \oplus
U_1)|_{\clq_P} \in \clb(\clq_P),\]is the compressed scalar operator.
Finally, since $w(A) \leq 1$, the decompressed numerical radius of
$X$ is not greater than one.

\NI For the converse part, let \[S = X + P X^*,\] for some
compressed scalar operator $X = \bm{P}_{\clq_P}(A \oplus
U_1)|_{\clq_P} \in \clb(\clq_P)$ with $w(A) \leq 1$. Then we
calculate
\[\begin{split}
S^* - S P^* & = (X + P X^*)^* - (X + P X^*) P^*\\ & = X^* + X P^* -
(X P^* + P X^* P^*) \\& = X^* - P X^* P^*.
\end{split}\]On the other hand, since \[ \begin{split}A
\oplus U_1 & = (I_{H^2(\mathbb{D})} \otimes A) \oplus U_1 \\& = ((\bm{P}_{\mathbb{C}} + M_z M_z^*) \otimes A) \oplus U_1\\
& = ((\bm{P}_{\mathbb{C}} \otimes A) \oplus 0) + ((M_z \otimes
I_{\cld_{P^*}}) \oplus U) ((I_{H^2(\mathbb{D})} \otimes A) \oplus
U_1) ((M_z \otimes I_{\cld_{P^*}}) \oplus U)^*,\end{split}\] and
\[\begin{split}& \bm{P}_{\clq_P}[((M_z \otimes I_{\cld_{P^*}}) \oplus U)
((I_{H^2(\mathbb{D})} \otimes A) \oplus U_1) ((M_z \otimes
I_{\cld_{P^*}}) \oplus U)^*]|_{\clq_P}\\ & = \bm{P}_{\clq_P}((M_z
\otimes I_{\cld_{P^*}}) \oplus U)(\bm{P}_{\clq_P^{\perp}} +
\bm{P}_{\clq_P}) ((I_{H^2(\mathbb{D})} \otimes A) \oplus
U_1)\bm{P}_{\clq_P} ((M_z \otimes I_{\cld_{P^*}}) \oplus
U)^*]|_{\clq_P}\\ & = \bm{P}_{\clq_P}((M_z \otimes I_{\cld_{P^*}})
\oplus U) \bm{P}_{\clq_P} ((I_{H^2(\mathbb{D})} \otimes A) \oplus
U_1)\bm{P}_{\clq_P} ((M_z \otimes I_{\cld_{P^*}}) \oplus
U)^*]|_{\clq_P}
\\ & = [ \bm{P}_{\clq_P}((M_z \otimes I_{\cld_{P^*}}) \oplus
U)]|_{\clq_P} [\bm{P}_{\clq_P} ((I_{H^2(\mathbb{D})} \otimes A)
\oplus U_1)|_{\clq_P}] [\bm{P}_{\clq_P} ((M_z \otimes
I_{\cld_{P^*}}) \oplus U)^*]|_{\clq_P} \\& = P [\bm{P}_{\clq_P} (A
\oplus U_1)|_{\clq_P}] P^*,
\end{split}\]
we have
\[\begin{split}\bm{P}_{\clq_P} (A
\oplus U_1)|_{\clq_P} & = \bm{P}_{\clq_P} ((I_{H^2(\mathbb{D})}
\otimes A) \oplus U_1)|_{\clq_P}\\ & = \bm{P}_{\clq_P}
(\bm{P}_{\mathbb{C}} \otimes A)|_{\clq_P} + P [\bm{P}_{\clq_P} (A
\oplus U_1)\bm{P}_{\clq_P}] P^*.
\end{split}\]
The last equality shows that \[X = \bm{P}_{\clq_P}
(\bm{P}_{\mathbb{C}} \otimes A)|_{\clq_P} + P X P^*,\]and hence
\[X^* - P X^* P^* = \bm{P}_{\clq_P} (\bm{P}_{\mathbb{C}} \otimes
A^*)|_{\clq_P}.\]It now follows that\[S^* - S P^* = \bm{P}_{\clq_P}
(\bm{P}_{\mathbb{C}} \otimes A^*)|_{\clq_P}.\] In view of Theorems
\ref{Aequiv} and \ref{A-equation}, we conclude that $(S, P)$ is a
$\Gamma$-contraction. This completes the proof of the theorem. \qed

\vspace{0.1in}

Let $(S, P)$ be a $\Gamma$-contraction. Then by the Wold
decomposition theorem of $\Gamma$-contractions (see the paragraph
after Theorem \ref{gamma-isometry} or Theorem 2.8 in \cite{AY5}) the
tuple can be decomposed as the direct sum of a c.n.u.
$\Gamma$-contraction and a $\Gamma$-unitary on the Wold
decomposition space of the contraction $P$. Consequently, the above
characterization is valid for any $\Gamma$-contraction on the Wold
decomposition space of the contraction where the operator $X$ is the
direct sum of a compressed scalar operator with the decompressed
numerical radius not greater than one with a unitary operator. For
$C_{.0}$ case, the necessary part of the above result was obtained
in \cite{BPR} in a different point of view.

Our next theorem concerns the uniqueness of the compressed scalar
operators of the representations of $\Gamma$-contractions in Theorem
\ref{XP}.

\begin{Theorem}\label{geom-2}
Let $(S, P)$ be a c.n.u. $\Gamma$-contraction on $\clq_P$. Then
there exists a unique compressed scalar operator $X$ with
$w_{\clq_P}(X) \leq 1$ such that \[S = X + P X^*.\]
\end{Theorem}

\NI\textsf{Proof.} The representation $S = X + P X^*$ for some
compressed scalar $X = \bm{P}_{\clq_P} (A \oplus U_1)|_{\clq_P}$
follows
from Theorem \ref{XP}. For the remaining part, we calculate \[\begin{split}S & = X + P X^*\\
& = \bm{P}_{\clq_P} (A \oplus U_1)|_{\clq_P} +  \bm{P}_{\clq_P} (M_z
\oplus U) \bm{P}_{\clq_P} (A \oplus U_1)^*|_{\clq_P}\\ & =
\bm{P}_{\clq_P} (A \oplus U_1)|_{\clq_P} +  \bm{P}_{\clq_P} (M_z
\oplus U) (I - P_{{\clq_P}^{\perp}}) (A \oplus U_1)^*|_{\clq_P} \\ &
= \bm{P}_{\clq_P} (A \oplus U_1)|_{\clq_P} +  \bm{P}_{\clq_P} (M_z
\oplus U) (A \oplus U_1)^*|_{\clq_P} \\ & = \bm{P}_{\clq_P} [(A
\oplus U_1) + (A \oplus U_1)^*(M_z \oplus U)]|_{\clq_P} \\ & =
\bm{P}_{\clq_P} [(A + A^* M_z) \oplus (U_1 + U_1^* U)]|_{\clq_P},
\end{split}\]where $U = M_{e^{it}}|_{\overline{\Delta_P
L^2_{\cld_P}(\mathbb{T})}}$. Consequently,  \[((A + A^* M_z) \oplus
(U_1 + U_1^* U), M_z \oplus U)\]is a $\Gamma$-isometric dilation of
$(S, P)$ and hence the result follows from the uniqueness part of
the $\Gamma$-isometric dilation of Theorem \ref{lift}. \qed

\newsection{Concluding Remarks}

\NI\textsf{(I) $\Gamma$-contractive Hilbert modules:}

Let $\{T_1, T_2\}$ be a pair of commuting operators in $\clb(\clh)$.
Then $\clh$ is a \textit{Hilbert module} (see \cite{DP}) over
$\mathbb{C}[z_1, z_2]$ where
\[ p \cdot h = p (T_1, T_2)h,\]for all $p \in
\mathbb{C}[z_1, z_2]$ and $h \in \clh$. A Hilbert module $\clh$ over
$\mathbb{C}[z_1, z_2]$ is said to be \textit{$\Gamma$-contractive
Hilbert module} if the ordered pair $(T_1, T_2)$ is a
$\Gamma$-contraction.

Examples of $\Gamma$-contractive Hilbert modules are :

(i) $\Gamma$-isometries,

(ii) $\Gamma$-unitaries,

(iii) direct sum of (i) and (ii).

Finally, by Theorem \ref{gamma-dilation-AY},

(iv) a $\Gamma$-contractive Hilbert module can be realized as a
compression of any one of (i), (ii) or (iii) to a joint co-invariant
subspace.

Now we turn to the class of pure $\Gamma$-isometries. Given a pure
$\Gamma$-isometry $(A + A^* M_z, M_z)$ on $H^2_{\cle_*}(\mathbb{D})$
for some $A \in \clb(\cle_*)$ with $w(A) \leq 1$, we say that
$H^2_{\cle_*}(\mathbb{D})$ is a \textit{$\Gamma$-isometric Hardy
module} with symbol $A$, where
\[p \cdot h = p(A + A^* M_z, M_z) h,\]for all $p \in \mathbb{C}[z_1,
z_2]$ and $h \in \clh$. Let $H^2_{\cle_*}(\mathbb{D})$ be a
$\Gamma$-isometric Hardy module with symbol $A$ and $\cls$ be a
closed subspace of $H^2_{\cle_*}(\mathbb{D})$. Then $\cls$ is said
to be a submodule of the $\Gamma$-isometric Hardy module
$H^2_{\cle_*}(\mathbb{D})$ if $\cls$ is invariant under $A + A^*
M_z$ and $M_z$.

Here we will present our Beurling-Lax-Halmos type theorem (Theorem
\ref{BLH}) in the Hilbert modules language.

\begin{Theorem}\label{BLH-mod}
Let $\cls \neq \{0\}$ be a closed subspace of
$H^2_{\cle_*}(\mathbb{D})$. Then $\cls$ is a submodule of the
$\Gamma$-isometric Hardy module $H^2_{\cle_*}(\mathbb{D})$ with
symbol $A$ if and only if there exists a $\Gamma$-isometric Hardy
module $H^2_{\cle}(\mathbb{D})$ with a symbol $B$ on $\cle$ and an
isometric module map
\[U : H^2_{\cle}(\mathbb{D}) \longrightarrow
H^2_{\cle_*}(\mathbb{D}),\]such that $\cls = U
H^2_{\cle}(\mathbb{D})$. Moreover, when such a $\Gamma$-isometric
Hardy module exists it is unique (up to unitary equivalence).
\end{Theorem}

One consequence of the Beurling-Lax-Halmos theorem is that a
non-zero submodule of the Hardy module $H^2_{\cle_*}(\mathbb{D})$ is
unitarily equivalent to $H^2_{\cle}(\mathbb{D})$ for some Hilbert
space $\cle$. This phenomenon is no longer true in general when one
consider the Hardy modules over the unit ball or the unit polydisc
in $\mathbb{C}^n$, $n \geq 2$. However, a non-zero submodule of a
$\Gamma$-isometric Hardy module is unitarily equivalent to a
$\Gamma$-isometric Hardy module.

\begin{Corollary}
A non-zero submodule of a $\Gamma$-isometric Hardy module is
isometrically isomorphic with a $\Gamma$-isometric Hardy module.
\end{Corollary}

Let $\cls = M_{\Theta}H^2_{\cle}(\mathbb{D})$ be a $M_z$-invariant
subspace of $H^2_{\cle_*}(\mathbb{D})$ for some inner multiplier
$\Theta \in H^{\infty}_{\clb(\cle, \cle_*)}(\mathbb{D})$ and that
$\cls$ be invariant under the multiplication operator $M_{p}$ where
$p$ is a $\clb(\cle_*, \cle)$-valued analytic polynomial. Then
\[ p \,\Theta = \Theta \Psi,\]for some unique $\Psi \in H^{\infty}_{\clb(\cle_*,
\cle)}(\mathbb{D})$.

\NI \textsf{Problem:} What is the representation of the unique
multiplier $\Psi$? Under what conditions that $\Psi$ will be a
polynomial, or a polynomial of the same degree of $p$?
\vspace{0.1in}

Theorem \ref{BLH-mult} implies that the question has a complete
answer when
\[p(z) = A + A^* z.\]One possible approach to solve this problem is
to consider first the finite dimension case, that is, $\cle_* =
\mathbb{C}^m$ for $m >1$.

Also one can formulate the above problem in the Hilbert modules
point of view. In this case, an isometric module map may yield a
natural candidate for $\Psi$. At present, we do not have any
positive result along that line.

\vspace{0.1in}

\NI\textsf{(II) Complete unitary invariants:}

We now turn to the task of determining a complete set of unitary
invariants of the class of c.n.u. $\Gamma$-contraction.

Let  $(S_1, P_1)$ on $\clh_1$ and $(S_2, P_2)$ on $\clh_2$ be a pair
of $\Gamma$-contractions. We consider the representation of $(S_i,
P_i)$ in $\clq_{P_i}$. Then by Theorem \ref{geom-2}, there exists
unique compressed scalar operators $X_i \in \clb(\clq_{P_i})$ with
$w_{\clq_{P_i}}(X)$, $i = 1, 2$, such that
\[S_i = X_i + P_i X_i^*. \quad \quad (i=1, 2)\]Furthermore, recall
that (see the proof of Theorem \ref{XP}) \[S_i^* - S_i P_i^* = X_i^*
= P_{\clq_{P_i}} (A_i^* \oplus U_i^*)|_{\clq_{P_i}} \in
\clb(\clq_P).\]We are now ready to prove that $(X, P)$ is a complete
set of unitary invariants for the class of c.n.u. contractions $(S,
P) = (X + P X^*, P)$, realized in the model space $\clq_P$.

\begin{Theorem}
$(X, P)$ is a complete set of unitary invariants for the class of
c.n.u. $\Gamma$ contractions. That is, for $(S_i, P_i)$ on
$\clq_{P_i}$, a pair of c.n.u. $\Gamma$-contractions on the model
space $\clq_{P_i}$, $(S_1, P_1) \cong (S_2, P_2)$ if and only if
\[(X_1, P_1) \cong (X_2, P_2).\]
\end{Theorem}
\NI \textsf{Proof.} Let $(X_1, P_1) \cong (X_2, P_2)$, that is,
\[X_2 = \bm{\tau} X_1 \bm{\tau}^*, \quad \mbox{and} \quad P_2 = \bm{\tau}
P_1 \bm{\tau}^*,\]for some unitary operator $\bm{\tau} \in
\clb(\clq_{P_1}, \clq_{P_2})$. Then \[S_2 = X_2 + P_2 X_2^* =
\bm{\tau} (X_1 + P_1 X_1^*) \bm{\tau}^*,\]and hence \[(S_1, P_1)
\cong (S_2, P_2).\]Conversely, let \[S_2 = \bm{\eta} S_1 \bm{\eta}^*
\quad \mbox{and} \quad P_2 = \bm{\eta} P_1 \bm{\eta}^*,\]for some
unitary operator $\bm{\eta} \in \clb(\clq_{P_1}, \clq_{P_2})$. Then
\[X_2^* = S_2^* - S_2 P_2^* = \bm{\eta} (S_1^* - S_1 P_1^*)
\bm{\eta}^* = \bm{\eta} X_1^* \bm{\eta}^*,\]that is, $(X_1, P_1)
\cong (X_2, P_2)$. \qed

\vspace{0.1in}

\NI\textsf{(III) Solving the commutant lifting theorem:}

The commutant lifting theorem was first proved by D. Sarason
\cite{Sa} and then in complete generality by Sz.-Nagy and Foias (see
\cite{NF1} and \cite{NF}). Since then, it has been identified as one
of the most useful results in operator theory. Here we recall a
special case of the commutant lifting theorem. Let $P$ be a c.n.u.
contraction on $\clh$. Let $X$ commutes with $P$. First, we identify
$P$ with the compression of the multiplication operator on the
Sz.-Nagy and Foias model space $\clq_P$. Next, we consider the
representation of $X$ on $\clq_P$, which we again denote by $X$.
Then the commutant lifting theorem implies the following commutative
diagram

\setlength{\unitlength}{3mm}
\begin{center}
\begin{picture}(40,16)(0,0)
\put(15,3){$\clq_P$}\put(19,1.6){$X$} \put(19,10.8){$\tilde{X}$}
\put(23.3,3){$\clq_P$} \put(23, 10){$\clh_P$}\put(15, 10){$\clh_P$}
\put(23.8,9.2){ \vector(0,-1){5}}\put(15,9.2){ \vector(0,-1){5}}
\put(17.3, 3.4){\vector(1,0){5.1}} \put(17.3,
10.4){\vector(1,0){5.1}}\put(11.9,7.3){$\bm{P}_{\clq_P}$}\put(25.3,7){$\bm{P}_{\clq_P}$}
\end{picture}
\end{center}
where $\tilde{X}$ commutes with the multiplication operator on
$\clh_P$ and $\|X\| = \|\tilde{X}\|$. It is usually a difficult
problem to find a solution $\tilde{X}$ to a given $X$ in the
commutator of $P$. One way to explain one of our main results,
namely, Theorem \ref{lift} is that if $S$ commutes with a given
c.n.u. contraction and if $\Gamma$ is a spectral set of $(S, P)$
(that is, $(S, P)$ is a $\Gamma$-contraction) then the solution to
the commutant lifting theorem is unique and explicit.

Therefore the results of this paper along with the seminal work of
Agler and Young (\cite{AY1} - \cite{AY6}) is an evidence of solving
the commutant lifting theorem uniquely and explicitly for a class of
commutators of a contraction.

Another possible approach to obtain some of the results of this
paper is to develop an independent proof of the characterization
result, Theorem \ref{XP}. Here, however, we do not pursue this
direction. Also we believe that our methods will be applicable not
only to other studies, but also demonstrate one way to set up and
solve the commutant lifting theorem in a more general framework.

Finally, following the work of Sz.-Nagy and Foias and by virtue of
our results, one can develop a $H^{\infty}$-functional calculus on
$\{(z_1+z_2, z_1 z_2) : |z_1|, |z_2| < 1\}$ of a c.n.u.
$\Gamma$-contraction. Moreover, a study of invariant subspaces of
$\Gamma$-contractions can be carried out. This will be considered in
future work.

\vspace{0.2in}

\NI \textsf{Acknowledgement:} We are grateful to Tirthankar
Bhattacharyya and Haripada Sau for pointing out some errors in the
earlier versions of this paper. We are also grateful to the referee
for reading the paper carefully and providing a number of valuable
suggestions.

\end{document}